\documentclass{cocv}
\usepackage{latexcad}

\renewcommand{\l}{\left}
\renewcommand{\r}{\right}
\newcommand{\f}[2]{\frac{#1}{#2}}
\newcommand{\inv}[1]{\frac{1}{#1}}
\renewcommand{\it}{\textit}
\renewcommand{\rm}{\textrm}
\newcommand{\no}{\noindent}
\newcommand{\ds}{\displaystyle}
     
\begin{document}
\title{Asymptotics of accessibility sets along an
abnormal trajectory}
\author{E. Tr\'elat}
\address{Universit\'e de Bourgogne, Laboratoire de Topologie, \\
UMR 5584 du CNRS, BP47870, 21078 Dijon Cedex, France  \\
e-mail : trelat@topolog.u-bourgogne.fr}
%
%
\begin{abstract} 
We describe precisely, under generic conditions, the contact of
the accessibility set at time $T$ with an abnormal direction,
first for a single-input affine control system with constraint on
the control, and then as an
application for a sub-Riemannian system of rank 2. As a
consequence we obtain in sub-Riemannian geometry a new
splitting-up of the sphere near an abnormal minimizer $\gamma$
into two sectors, bordered by the first Pontryagin's cone along
$\gamma$, called the $\xLinfty$-sector and the $\xLtwo$-sector.
Moreover we find again necessary and sufficient conditions of
optimality of an abnormal trajectory for such systems,
for any optimization problem.
\end{abstract}
\begin{resume}
On d\'ecrit pr\'ecis\'ement, sous des conditions g\'en\'eriques,
le contact des ensembles accessibles en temps $T$ avec une
direction anormale, tout d'abord pour un syst\`eme de contr\^ole
affine mono-entr\'ee avec contrainte sur le contr\^ole, puis
comme application pour un syst\`eme sous-Riemannien de rang 2.
Comme cons\'equence on obtient en g\'eom\'etrie sous-Riemannienne
un nouveau d\'ecoupage de la sph\`ere au voisinage d'une anormale
minimisante $\gamma$ en deux secteurs appel\'es secteur
$\xLinfty$ et secteur $\xLtwo$, d\'elimit\'es par le premier c\^one
de Pontryagin le long de $\gamma$. De plus on retrouve des
conditions n\'ecessaires et suffisantes d'optimalit\'e pour une
trajectoire anormale de tels syst\`emes, pour un probl\`eme
d'optimisation quelconque.
\end{resume}
\subjclass{???, ???}
\maketitle


\section{Introduction}
\subsection{Abnormal trajectories}
Consider a control system on $\xR^n$~:
\begin{equation} \label{systgen1}
\dot{x}(t)=f(x(t),u(t)),\ 
x(0)=x_0
\end{equation}
where $f:\xR^n\times\xR^m\longrightarrow \xR^n$ is smooth,
$x_0\in\xR^n$,
and the set of admissible controls ${\cal U}$ is made of
measurable bounded functions
$u~:[0,T(u)]\longrightarrow\Omega\subset\xR^m$. 

\begin{dfntn}
Let $T>0$. The \it{end-point mapping} at time $T$ of system
(\ref{systgen1}) is the mapping
$$E_T~:\begin{array}{rcl}
{\cal{U}}&\longrightarrow&\xR^n\\
u&\longmapsto&x_u(T) \end{array}$$
where $x_u$ is the trajectory associated to $u$.
\end{dfntn}

It is a very classical fact that $E_T$ is smooth in the
$\xLinfty$ topology, if ${\cal U}\subset \xLinfty([0,T])$.

\begin{dfntn}
The \it{time$\times$end-point mapping} of system
(\ref{systgen1}) is the mapping
$$F~:\begin{array}{rcl}
\xR^+\times{\cal{U}}&\longrightarrow&\xR^n\\
(T,u)&\longmapsto&E_T(u) \end{array}$$
\end{dfntn}

\begin{dfntn} \label{defi1}
A control $u$ (or the corresponding trajectory $x_u$)
is said to be \it{abnormal} on $[0,T]$ if $(T,u)$ is a
singular point of the mapping $F$.
\end{dfntn}

An equivalent definition may be given using the well-known
\it{Pontryagin Maximum Principle} (see \cite{P}), which
parametrizes the previous condition~:

\begin{dfntn} \label{defi2}
A control $u$ on $[0,T]$ (or the corresponding trajectory $x_u$)
is said to be \it{abnormal} if
there exists a non trivial absolutely continuous
function $p(.)~: [0,T]\longrightarrow \xR^n\times\xR$
called adjoint vector,
such that the trajectory $x$ satifies almost everywhere the
system~:
\begin{equation} \label{systPMP}
\dot{x}=\f{\partial H}{\partial p},\ 
\dot{p}=-\f{\partial H}{\partial x},\ 
\f{\partial H}{\partial u}=0
\end{equation}
where $H(x,p,u)=<p,f(x,u)>$ is the Hamiltonian of
the system, and moreover~:
\begin{equation} \label{annulationHam}
H(x,p,u)=0 \quad \rm{a.e. on }[0,T]
\end{equation}
\end{dfntn}

\begin{rmrk}
If a control $u$ is abnormal on $[0,T]$ then it is abnormal on
$[0,t]$ for any $t\in[0,T]$.
\end{rmrk}

\begin{dfntn}
An abnormal control is said to be \it{of corank 1} if the adjoint
vector $p(.)$ is defined uniquely up to a scalar multiplier.
\end{dfntn}

\begin{dfntn}
Let $u$ be an abnormal control on $[0,T]$, and $x_u$ its
associated trajectory.
The subspace $\xim dE_t(u)$ is called the \it{first
Pontryagin's cone} at $x_u(t)$.
\end{dfntn}

An abnormal control is of corank 1 on $[0,T]$ if the subspace
$\xim dE_T(u)$ has codimension 1 in $\xR^n$.

\begin{rmrk}
An abnormal control on $[0,T]$ is in particular a singular point
of the end-point mapping $E_T$. Conversely if $u$ is a
singularity of the end-point mapping and if moreover the
trajectory $x_u$ associated to $u$ is almost everywhere tangent
to its first Pontryagin's cone then $u$ is abnormal.
\end{rmrk}


\subsection{Accessibility sets}
\begin{dfntn}
Consider the control system (\ref{systgen1}), and let
$T>0$. The \it{accessibility set at time $T$}, denoted by
$Acc(T)$, is the set of points that can be reached from $x_0$ in
time $T$ by solutions of system (\ref{systgen1}), i.e. this is
the image of the end-point mapping $E_T$.
\end{dfntn}

Let $\gamma$ be a reference trajectory on $[0,T]$, solution of
(\ref{systgen1}), associated to a control $u$. Our aim is to
describe $Acc(T)$ near $\gamma(T)$. If $u$ is not a singular
point of the end-point mapping, then
obviously $Acc(T)$ is open near $\gamma(T)$.
The situation when $u$ is a critical point of $E_T$ has to be
analyzed.


\subsection{Single-input affine systems}
Consider a smooth \it{single-input affine system in $\xR^n$ with
constraint on the control}~:
\begin{equation} \label{introsyst2}
\begin{split}
\dot{x}(t)&=X(x(t))+u(t)Y(x(t)),\quad |u(t)|\leq\eta\\
x(0)&=x_0
\end{split}
\end{equation}
where $X, Y$ are smooth vector fields in $\xR^n$, and $u$ is a
scalar measurable function.
Suppose that the trajectory $\gamma$ associated to the control
$u=0$ is abnormal and of corank 1
on $[0,T]$. Let $p(.)$ be an adjoint vector
associated to $\gamma$. Then the condition $H=0$ is equivalent
to~:
$$<p(t),X(\gamma(t))>=0 \ \rm{ a.e.}$$
i.e. $X\in\xim dE_t(u)$ along $\gamma$.

The aim of this article is to describe precisely the accessibility
set $Acc(T)$ at time $T$ near $\gamma(T)$.

The basic object we have to study
is the so-called \it{intrinsic second-order derivative}
of the end-point mapping~:

\begin{dfntn}
The \it{intrinsic second-order derivative} along $\gamma$ is the
real quadratic form~:
$$E_u''(v)=p(T).d^2E_T(u).(v,v)$$
where $v\in\xker dE_T(u)$.
\end{dfntn}

We shall make a spectral
analysis of the end-point mapping along such an abnormal
trajectory. Using the formalism and normal forms of \cite{BK} we
shall represent the intrinsic second-order derivative by an
explicit differential operator along the abnormal trajectory
$\gamma$~; then a spectral analysis of this
operator shall lead to a precise description of the boundary of
$Acc(T)$ near $\gamma(T)$ (Theorem \ref{bordaffine}).
In Section \ref{sectionsrcase} we apply this result to
sub-Riemannian systems of rank 2 (Theorem \ref{bordSR}),
and obtain a new splitting-up
of the sub-Riemannian sphere near an abnormal minimizer $\gamma$
into two sectors, in which the behaviours of minimizing
trajectories near $\gamma$ are topologically different (Section
\ref{secteursSR}). On the other part this theory on accessibility
sets leads to find
again some well-known results on optimality of abnormal
trajectories, that we recall and improve slightly in Section
\ref{sectionoptimalité}.


\section{Asymptotics of the accessibility sets}
In this Section we describe precisely the boundary of
accessibility sets for a single-input affine system with
constraint on the input near a reference abnormal trajectory.
Then we apply our results to the
sub-Riemannian case of rank 2 in order to get the
\it{contact} of the sphere with the abnormal direction. As a
consequence we obtain a \it{splitting-up of the sphere into
two sectors near the abnormal minimizer}.

\subsection{Single-input affine control systems}
Consider a smooth \it{single-input affine control system} in
$\xR^n$, $n\geq 3$~:
\begin{equation} \label{af}
\dot{x}(t)=X(x(t))+u(t)Y(x(t)),\ x(0)=0
\end{equation}
with the \it{constraint} on the control
\begin{equation} \label{contrainte}
|u(t)|\leq\eta
\end{equation}
Let $Acc^\eta(T)$ denote the \it{accessibility set} at time $T$
for this affine system with constraint $\eta$ on the control.
Let $\gamma$ be a reference trajectory defined on
$[0,T]$. In the sequel we make the following assumptions along
$\gamma$~:
\begin{itemize}
\item[($H_0$)] $\gamma$ is injective, associated to $u=0$ on
$[0,T]$.
\item[($H_1$)] $\forall t\in[0,T] \quad K(t)=\rm{Vect }\{
ad^kX.Y(\gamma(t))\   / \   k\in\xN\}$ (first
Pontryagin's cone along $\gamma$) has codimension 1, and is
spanned by the first $n-1$  vectors, i.e.~:
$$K(t)=\rm{Vect }\{ad^kX.Y(\gamma(t))\   / \   k=0\ldots n-2\}$$
\item[($H_2$)] $\forall t\in[0,T]\quad ad^2Y.X(\gamma(t))\notin
K(t)$.
\item[($H_3$)] $\forall t\in[0,T] \quad X(\gamma(t))\notin
\rm{Vect }\{ad^kX.Y(\gamma(t))\   / \   k=0\ldots n-3\}$.
\item[($H_4$)] $\forall t\in[0,T] \quad X(\gamma(t))\in K(t)$.
\end{itemize}
In these conditions $\gamma$ is \it{abnormal and of corank 1}.
Actually assumptions $(H_1-H_3)$ are generic,
see \cite{BK2}.
Moreover we get normal forms in which the intrinsic
second-order derivative may be represented by an explicit
differential operator (see \cite{BK}), in the following way. 

\begin{lmm}\cite{BK} \label{formenormaleexcept}
Under the previous assumptions, the system $(X,Y)$ is
in a $\xCzero$-neighborhood of $\gamma$
feedback-equivalent to~:
\begin{equation} \label{systemAN}
\begin{split}
f_0&=\f{\partial}{\partial x_1}+\sum_{i=1}^{n-2}x_{i+1}
\f{\partial}{\partial x_i} +\sum_{i,j=2}^na_{ij}(x_1)x_ix_j
\f{\partial}{\partial x_n} +\sum_{i=1}^{n-1}x_nf_i(x_1)
\f{\partial}{\partial x_i} +\sum_{i=1}^nZ_i \f{\partial}{\partial
x_i} \\
f_1&=\f{\partial}{\partial x_{n-1}}
\end{split}
\end{equation}
where $a_{n-1,n-1}(t)>0$ on $[0,T]$, and the 1-jet
(resp. 2-jet) of $Z_i,i=1\ldots n-1$ (resp. $Z_n$) along $\gamma$
is equal to $0$.

\no Moreover the feedback $(\varphi ,\alpha ,\beta )$ satisfies~:
$$(\varphi ,\alpha ,\beta ).(X,Y)=(f_0,f_1)$$
where $f_0=\varphi _*(X+\alpha Y)$, $f_1=\varphi _*(\beta Y)$,
and $\varphi $ is a germ of diffeomorphism along $\gamma$ such
that~:
\begin{itemize}
\item[(i)] $\varphi (x_1,0,\ldots,0)=(x_1,0,\ldots,0)$
\item[(ii)] $\f{\partial \varphi }{\partial
x_{n-1}}=(0,\ldots,0,*,0)$
\end{itemize}
and $\alpha ,\beta $ are real functions defined in a neighborhood
of $\gamma$ such that $\beta$ does not vanish along $\gamma$
and $\alpha _{|\gamma}=0$.
\end{lmm}

The k-jet is defined as follows~:

\begin{dfntn}
Let $\ds{V=\sum_{k=1}^nV_k\f{\partial}{\partial x_k} }$ a vector
field. Since $\gamma$ is given by $x_1(t)=t$ and $x_i(t)=0,
i=2\ldots n$, the component $V_k$ can be written in a
neighborhood of $\gamma$ as $\ds{\sum_{p=1}^{+\infty}j_pV_k}$,
where $j_0V=V_{/\gamma}$, and~:
$$j_1V_k=\sum_{i=2}^n a_i^k(x_1)x_i,\ 
j_2V_k=\sum_{i,j=2}^nb_{ij}^k(x_1)x_ix_j,\ \ldots$$
Set $\ds{j_iV=\sum_{k=1}^nj_iV_k\f{\partial}{\partial x_k} }$.
Then $\ds{\sum_{i=0}^kj_iV }$ is called \it{the k-jet of $V$
along $\gamma$}.
\end{dfntn}

Set $x_1=t+\xi$. In these conditions,
the controllable part of the system is
$(\xi,x_2,\ldots, x_{n-1})$, the reference abnormal trajectory is
$\gamma(t)=(t,0,\ldots,0)$, and the intrinsic second-order
derivative $E''_0(v)$
along $\gamma$ is identified to~:
\begin{equation*}
E''_0(v)=
\int_0^T\sum_{i,j=2}^{n-1} a_{ij}(t)\xi_i(t)\xi_j(t)\ dt\ , \rm{
where~:}
\end{equation*}
$$\dot{\xi}_1=\xi_2,\ldots,
\dot{\xi}_{n-2}=\xi_{n-1},\dot{\xi}_{n-1}=v,\quad\rm{and}\quad
\xi_i(0)=\xi_i(T)=0, i=1\ldots,n-1.$$
Integrating by parts, it can be
written into two different ways, namely either as a quadratic
form in $\xi_1$ or as a quadratic form in $\xi_2$~:

\begin{enumerate}
\item It is equal to ${Q_1}_{/G_1}(\xi_1)$, where~:
$$Q_1(\xi_1)=\int_0^Tq_1(\xi_1)dt\quad\rm{and}\quad
q_1(\xi_1)=\sum_{i,j=1}^{n-2}b_{ij}\xi_1^{(i)}\xi_1^{(j)}$$
with $b_{i-1,j-1}=\f{a_{ij}+a_{ji}}{2}$, and where $G_1$ is
the following space corresponding to the kernel
$\xker dE_T(0)$ of the first
derivative of the end-point mapping~:
$$G_1=\{\xi_1\ /\ \xi_1^{(2(n-2))}\in \xLtwo([0,T]),\ 
\xi_1^{(i)}(0)=\xi_1^{(i)}(T)=0,\  i=0\ldots n-2 \}$$

\begin{lmm} \label{lem1}
The quadratic form $Q_1$ is represented on $G_1$ by the operator
$D_1$ so that~:
$$Q_1(\xi_1)=(\xi,D_1\xi)_{\xLtwo}$$
where $(\ ,\ )$ is the usual scalar product in $\xLtwo([0,T])$,
and~:
\begin{equation} \label{opD}
D_1\ =\ \inv{2}\sum_{i=1}^{n-2}(-1)^i\f{d^i}{dt^i}\f{\partial
q_1}{\partial y^{(i)}}
\ =\ \sum_{i,j=1}^{n-2}(-1)^j\f{d^j}{dt^j}b_{ij}\f{d^i}{dt^i}
\end{equation}
\end{lmm}

\item It is equal to ${Q_2}_{/G_2}(\xi_2)$, where~:
$$Q_2(\xi_2)=\int_0^Tq_2(\xi_2)dt\quad\rm{and}\quad
q_2(\xi_2)=\sum_{i,j=0}^{n-3}b_{i+1,j+1}x_2^{(i)}x_2^{(j)}$$
and where $G_2$ is the space corresponding to the kernel
of the first derivative~:
\begin{equation*}
G_2=\{\xi_2\ /\ \xi_2^{(2(n-3))}\in \xLtwo([0,T]),\ 
\xi_2^{(i)}(0)=\xi_2^{(i)}(T)=0,\ i=0\ldots n-3, \rm{ and }
\int_0^T\xi_2\ dt=0 \}
\end{equation*}

\begin{lmm} \label{lem2}
The quadratic form $Q_2$ is represented on $G_2$ by the operator
$D_2$ so that~:
$$Q_2(\xi_2)=(\xi_2,D_2\xi_2)_{\xLtwo}$$
where
\begin{equation} \label{opD1}
D_2\ =\ \inv{2}\sum_{i=0}^{n-3}(-1)^i\f{d^i}{dt^i}\f{\partial
q_2}{\partial y^{(i)}}
\ =\ \sum_{i,j=0}^{n-3}(-1)^j\f{d^j}{dt^j}b_{i+1,j+1}\f{d^i}{dt^i}
\end{equation}
\end{lmm}
\end{enumerate}

Note that $Q_1(\xi)=Q_2(\dot{\xi})$ and $\ds{D_1=-\f{d}{dt}
D_2 \f{d}{dt} }$.

Our aim is to make a spectral analysis of these operators $D_1,
D_2$. Unfortunately the spectrum of $D_i$ on $G_i$ is empty.
Hence we shall enlarge the Sobolev space $G_i$ so that the
spectrum is not trivial and that Representation Lemmas \ref{lem1}
and \ref{lem2} are still valid. We set~:

\begin{itemize}
\item $F_1=\{\xi_1\ /\ \xi_1^{(n-2)}\in \xLtwo([0,T]),\ 
\xi_1^{(i)}(0)=\xi_1^{(i)}(T)=0,\  i=0\ldots n-3 \}$ for the
operator $D_1$. Endowed with the norm $||\xi_1||_{F_1}
=||\xi_1^{(n-2)}||_{\xLtwo}$, $F_1$ is a Sobolev space.
\item $F_2=\{\xi_2\ /\ \xi_2^{(n-3)}\in \xLtwo([0,T]),\ 
\xi_2^{(i)}(0)=\xi_2^{(i)}(T)=0,\  i=0\ldots n-4 \}$ for $D_2$ if
$n\geq 4$ (if $n=3$, no condition is imposed).
Endowed with the norm $||\xi_2||_{F_2}
=||\xi_2^{(n-3)}||_{\xLtwo}$, $F_2$ is a Sobolev space.
\end{itemize}

\begin{dfntn}
Let $i=1$ or $2$.
We call $T$ a \it{conjugate time} of $Q_i$ along $\gamma$ is
there exists $\xi_i\in F_i\backslash \{ 0\}$ such that
$\xi_i^{(2(n-1-i))}\in \xLtwo([0,T])$ and $D_i\xi_i=0$.
\end{dfntn}

\begin{lmm}(see for instance \cite{Brezis})
For any $f\in \xLtwo([0,T])$, if $T$ is not a conjugate time,
there exists $\xi_i\in F_i$ unique such that
$\xi_i^{(2(n-1-i))}\in \xLtwo([0,T])$ and $D_i\xi_i=f$. Let $L$
denote the 
operator $f\mapsto \xi_i$ \it{considered as an operator from
$\xLtwo([0,T])$ into $\xLtwo([0,T])$}~; it is selfadjoint
and compact.
\end{lmm}

The following Lemma is an improvement of
\cite{BK}, where only a non strict inequality is proved~:

\begin{lmm} \label{tempsconj}
Let $t_c$ (resp. $t_{cc}$) denote the first conjugate time of $Q_1$
on $F_1$ (resp. $Q_2$ on $F_2$). We have~: $0<t_{cc}<t_c$.
\end{lmm}

\begin{proof}
This proof can be found in \cite{BK}, where only a non strict
inequality is proved.
It is included only for convenience of the reader.

Let $y\in F_1$ such that $y^{(2(n-2))}\in \xLtwo([0,T])$, and
$z=\dot{y}$~;
then $Q_1(y)=Q_2(z)$. Let $\lambda$ (resp. $\mu$) denote the
smallest eigenvalue of $D_1$ (resp. $D_2$).
From spectral properties of compact
selfadjoint operators, we have~:
$$\lambda=\inf_{y\in F_1} \f{Q_1(y)}{(y,y)}\quad\rm{ and }\quad
\mu=\inf_{z\in F_2} \f{Q_2(z)}{(z,z)}$$
By definition~:
$$\f{Q_1(y)}{(y,y)}=\f{Q_2(z)}{\int_0^T(\int_0^tz(s)ds)^2dt}$$
where $z\in F_2$ and $\int_0^Tz(t)dt=0$. Hence~:
$$\lambda\geq \inf\l\{ \f{Q_2(z)}{\int_0^T(\int_0^tz(s)ds)^2dt}
\ /\ z\in F_2\backslash\{0\}\r\}$$
Moreover we get from the Cauchy-Schwarz inequality~:
$$\int_0^T\l(\int_0^tz(s)ds\r)^2dt\leq \int_0^Tt\int_0^tz(s)^2ds$$
and integrating by parts~:
$$\int_0^Tt\int_0^tz(s)^2dsdt\leq \f{T^2}{2}\int_0^Tz(t)^2dt
-\inv{2}\int_0^T t^2z(t)^2<\f{T^2}{2}\int_0^Tz(t)^2dt$$
Therefore~:
$$\lambda>\f{2}{T^2}\inf\l\{\f{Q_2(z)}{\int_0^Tz(t)^2dt}\ /\ z\in
F_1\r\}=\f{2}{T^2}\mu$$
From \cite{Sar} we know that eigenvalues of $D_1$ (or $D_2$) are
continuous and decreasing functions of $T$.
Hence if $T=t_c$ then $\lambda=0$, thus $\mu<0$ and consequently
$t_{cc}<t_c$.
\end{proof}

\begin{rmrk} \label{remtcc}
If $n=3$, we have $t_{cc}=+\infty$ provided assumptions
($H_0-H_3$) are fulfilled on $\xR^+$.
\end{rmrk}

\begin{rmrk}
The notion of conjugate
time \it{does not depend on the
constraint} on the control. It comes from the fact that the
abnormal
reference control belongs to the \it{interior} of the domain of
constraints.
\end{rmrk}

The main result is the following.

\begin{thrm} \label{bordaffine}
Consider the affine system (\ref{af}) with the constraint
(\ref{contrainte}), and suppose that assumptions ($H_0-H_4$) are
fulfilled along the reference abnormal trajectory $\gamma$ on
$[0,T]$. Let
$t_{cc}$ and $t_c$ denote the first conjugate times associated to
$\gamma$. Then~:
\begin{enumerate}
\item There exist coordinates $(x_1,\ldots,x_n)$ locally along
$\gamma$ such that in these coordinates~:
$\gamma(t)=(t,0,\ldots,0)$, and the first Pontryagin's cone along
$\gamma $ is~: $K(t)=\rm{Vect
}\{\f{\partial}{\partial x_1},\ldots,\f{\partial}{\partial
x_{n-1}} \}_{|\gamma}$.
\item If $T$ is small enough then for any point $(x_1,\ldots,x_n)$
of $Acc^\eta(T)\backslash \{\gamma(T)\}$
close to $\gamma(T)$ we have~: $x_n>0$
(see Fig. \ref{figparab}).

\drawingscale{0.4mm}
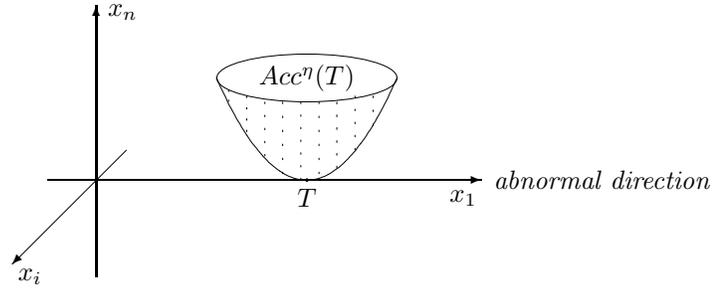
\placedrawing[h]{fparab.lp}{Shape of $Acc^\eta(T)$,
$T$ small}{figparab}

\item If $T<t_c$, then in the plane $(x_1,x_n)$, near the
point $(T,0)$, the boundary of $Acc^\eta(T)$ does not depend on
$\eta$, is a curve of class $\xCtwo$ tangent to the abnormal
direction, and its first term is~:
$$x_n=A_T(x_1-T)^2+\rm{o}((x_1-T)^2)$$
The function $T\mapsto A_T$ is continuous and strictly decreasing
on $[0,t_c[$. It is positive on $[0,t_{cc}[$ and negative on
$]t_{cc},t_c[$.\\
Moreover, if $\eta$ depends on $x_1-T$ then the result is still
valid providing~: $x_1-T=\rm{o}(\eta)$ as $x_1\to T$.

\item If $T>t_c$ then $Acc^\eta(T)$ is open near $\gamma(T)$.
\end{enumerate}
\end{thrm}

The evolution in function of $T$ of the intersection of
$Acc^\eta(T)$ with the
plane $(x_1,x_n)$ is represented on
Fig. \ref{figaccexc}. The contact with the abnormal direction is
of order 2~; the coefficient $A_T$ describes the \it{concavity}
of the curve. Beyond $t_c$ the accessibility set is open.

\drawingscale{0.5mm}
\placedrawing[h]{faccexc.lp}{}{figaccexc}

\begin{rmrk} \label{remcalculA_T}
The coefficient $A_T$ can be computed in the following way
(see \cite{BK}, and the Proof just below).
Actually there exists a function $J$ of class $\xCn{2(n-2)}$ on
$[0,T]$ such that $D_1J=0$ and satisfying the limit conditions~:
$$\forall k\in \{0,\ldots,n-3\}\quad J^{(k)}(0)=0,\
J^{(k)}(T)=\delta_{0}^k$$
Then~:
\begin{equation} \label{calculA_T}
A_T=Q_1(J)
\end{equation}
\end{rmrk}

\begin{rmrk}
Let us give a geometric interpretation of the role of the
quadratic forms $Q_1, Q_2$.
\begin{enumerate}
\item Role of $Q_1$. On the one part note
that $G_1$ is dense in $F_1$ and $Q_1$ is
continuous on $F_1$ for the norm $||\ ||_{F_1}$~;
hence the sign of $Q_1$ on $G_1$ is the same as on $F_1$. On the
other part, from the definition of the first conjugate time $t_c$
and the extremal properties of selfadjoint compact operators, we
get that ${Q_1}_{/F_1}$ is positive definite if $T<t_c$, and
indefinite if
$T>t_c$. Hence the same goes for ${Q_1}_{/G_1}$.

Now the
interpretation of ${Q_1}_{/G_1}$ is the following~: it is an
equivalent of the coordinate $x_n(T)$ as all others coordinates
are fixed~: $x_1(T)=T, x_2(T)=\cdots=x_{n-1}(T)=0$.

As a consequence if $T<t_c$ then points
$(T,0,\ldots,0,\varepsilon)$, where $\varepsilon\geq 0$ are
attainable at time $T$. If $T>t_c$, then the same goes with
$\varepsilon\geq 0$ or $\varepsilon\leq 0$.

\item Role of $Q_2$. Contrarily to the previous case, the
subspace $G_2$ is a strict subspace of $F_2$. However, set~:
$$H_2=\{\xi_2\ /\ {\xi_2}^{(2(n-3))}\in \xLtwo([0,T]),
{\xi_2}^{(i)}(0)={\xi_2}^{(i)}(T)=0, i=0\ldots n-3\}$$
Then $H_2$ is dense in $F_2$, and $Q_2$ is continuous on $F_2$ for
the norm $||\ ||_{F_2}$.

Hence using the same reasoning we obtain the following~: if
$T<t_{cc}$ then ${Q_2}_{/H_2}$ is positive definite~; if
$T>t_{cc}$ then ${Q_2}_{/H_2}$ is indefinite.

Now note that ${Q_2}_{/H_2}$ is an
equivalent of the coordinate $x_n(T)$ as coordinates
$x_2(T),\ldots,x_{n-1}(T)$ are fixed to $0$, but \it{the
coordinate $x_1(T)$ is let free}.

As a consequence if $T<t_{cc}$ then points
$(T\pm\delta,0,\ldots,0,\varepsilon)$, where $\varepsilon, \delta
\geq 0$ are
attainable at time $T$. If $T>t_{cc}$, then the same goes with
$\varepsilon\geq 0$ or $\varepsilon\leq 0$.

\end{enumerate}

This gives us the qualitative shape of $Acc(T)$. Actually this
interpretation is enough to study the time-optimality
of the abnormal trajectory $\gamma$, see Section
\ref{sectiontimeoptimality}.
Here Theorem \ref{bordaffine} gives more~: it describes
the contact of $Acc(T)$ with the abnormal direction.
\end{rmrk}

\begin{proof}[Proof of Theorem \ref{bordaffine}]
We shall recall briefly the formalism used by \cite{BK}, which
leads actually to a very precise decomposition of the intrinsic
second-order derivative. It is based on normal forms of Lemma
\ref{formenormaleexcept}.

Let $D_1^t$ denote the operator (\ref{opD}) and $D^t_2$ the
operator (\ref{opD1}).

\begin{lmm} \label{lemma0}
For all $t\in [0,T]$ there exists an Hilbertian basis $(e_n^t)
_{n\in\xN}$ of $\xLtwo([0,T])$ such that
\begin{itemize}
\item $\forall n\in\xN\quad e_n^t$ is $\xCn{2(n-2)}$ and $e_n^t
\in F_1$.
\item $\forall n\in\xN\quad D^t_1e_n^t=\lambda_n^te_n^t$.
\item $\lambda_1^t\leq\lambda_2^t\leq\cdots\leq\lambda_n^t\leq\cdots$
\end{itemize}
The eigenvalues $\lambda_n^t$ are continuous and decreasing functions
of $t$, and $\lambda_n^t\underset{n\rightarrow +\infty}
{\longrightarrow} +\infty$. Moreover, if $t_c$ denotes the first
conjugate time of $D^t_1$, then~:
\begin{itemize}
\item if $0\leq t< t_c$ then $\lambda_1^t>0$,
\item if $t>t_c$ then $\lambda_1^t<0$.
\end{itemize}
\end{lmm}

\begin{lmm} \label{lemma1}
Let $0<t<t_c$. Then there exist $J_i^t, \bar{J}_i^t$ in
$\xCn{2(n-2)}([0,t])$, $i=1\ldots n-2$, uniquely defined by the
following equations~:
\begin{itemize}
\item[(i)] $D^t_1J_i^t=D^t_1\bar{J}_i^t=0,\ i=1\ldots n-2$
\item[(ii)] $\begin{array}{l}
J_i^{(k)}(0)=\bar{J}_i^{(k)}(t)=0 \\
 J_i^{(k)}(t)=\bar{J}_i^{(k)}(0)=\delta_{i-1}^k,\ k=0\ldots n-3
\end{array}$
\end{itemize}
Every $\xi$ in $\xCn{2(n-2)}([0,t])$ can be expanded in a uniformly
convergent series~:
$$\xi=\sum_{i=1}^{n-2}\alpha_i{J}_i^t+\sum_{i=1}^{n-2}
\beta_i\bar{J}_i^t
+\sum_{n=1}^\infty y_ne_n^t$$
where $y_n\in\xR, \alpha_i=\xi^{(i-1)}(0), \beta_i=
\xi^{(i-1)}(t)$.
\end{lmm}

We shall use this lemma in the following way. \it{In first
approximation} the system written in the normal form of Lemma
\ref{formenormaleexcept} is given by~:
$$\dot{\xi}=x_2,\ \dot{x}_2=x_3,\ \ldots,\ \dot{x}_{n-1}=u$$
except the last coordinate, which is given at order 2 by~:
$$\dot{x}_n=\sum_{i,j=2}^{n-1}a_{ij}x_ix_j$$
Now the meaning of the $J_i, \bar{J}_i$'s is clear~: the functions
$\bar{J}_i$ (resp. $J_i$) represent the initial (resp. final)
conditions of the $n-2$ first coordinates. More precisely we have,
see \cite{BK}~:

\begin{lmm} \label{lemma2}
Let $T<t_{c}$, and $\ds{\xi=\sum_{i=1}^\infty v_ie_i^T
+\sum_{i=1}^{n-2}\alpha_iJ_i^T}$. We have~:
\begin{itemize}
\item
$\ds{E_i^T(\xi^{(n-1)})=\delta_1^iT+\alpha_i
+\rm{o}(||\xi^{(n-1)}||_{\xLinfty})}$, $i=1\ldots n-2$.
\item
$\ds{E_n^T(\xi^{(n-1)})=Q_1^T(\xi)+\rm{o}(||\xi^{(n-1)})||^2
_{\xLinfty})}$, where
$$Q_1^T(\xi)=
\sum_{i,j=1}^{n-2}A_{ij}^T\alpha_i\alpha_j
+\sum_{i=1}^\infty \lambda_i^Tv_i^2$$
and $A_{ij}^T=\check{Q}_1^T(J_i^T,J_j^T)$, where $\check{Q}^T_1$
denotes the symmetric bilinear form associated to the quadratic
form $Q^T_1$. In fact~:
$\ds{Q_1^T(\xi)=d^2E_n^T(0).(u,u)}$, where $u\sim\xi^{(n-1)}$, i.e.
\begin{equation} \label{onestar}
Q_1^T(\xi)=\int_0^T\sum_{i,j=1}^{n-2}b_{ij}(t)\xi^{(i)}(t)
\xi^{(j)}(t)dt
\end{equation}
where $b_{n-2,n-2}$ is continuous and strictly positive on $[0,T]$.
\end{itemize}
\end{lmm}

Hence with this formalism we study the image of the
end-point 
mapping by constructing directly the trajectories (and not by
choosing 
controls). Actually once a function $\xi$ has been determined,
the corresponding control is
$u=\xi^{(n-1)}+\rm{o}(||\xi^{(n-1)}||_{\xLinfty})$.

First of all using expression (\ref{onestar}) it is easy to prove
point 2 of Theorem \ref{bordaffine}. Indeed to study the
accessibility set at time $T$ from $0$ we have to consider
functions $\xi$ such that~:
$$\xi^{(i)}(0)=0,\ i=1\ldots n-2$$
Hence using repeatedly Poincar\'e's inequality\footnote{$f(0)=0
\implies ||f||_{\xLinfty([0,T])}\leq \sqrt{T}
||f'||_{\xLtwo([0,T])}$} we check that \it{if $T$ is small enough}
then there exists $\beta>0$ such that~:
$$Q_1^T(\xi)\geq \beta\int_0^T (\xi^{(n-2)}(t))^2dt$$
We shall deduce the result by proving that actually $E_n^T(u)>0$
for the corresponding control $u$. The result would be immediate
if $E_n^T$ were $\xCtwo$ in $H^{-1}$, i.e. in $\xLtwo$ topology
for coordinate $x_{n-1}$. But this is wrong (for instance the
function $f\mapsto \int_0^T(f^2+f^3)dt$ is not $\xCtwo$ on $\xLtwo$).
Anyway $E_n^T$ is a little better than $\xCone$. Precisely we
have (see \cite{AS}, Theorem 7.1)~:

\begin{lmm}
$E_n^T(u)=Q_1^T(\xi)+R(\xi)$ where
$|R(\xi)|\leq C||\xi^{(n-2)}||_{\xLinfty}||\xi^{(n-2)}||^2_{\xLtwo}$.
\end{lmm}

And hence if $T$ is small enough then
for any point $(x_1,\ldots,x_n)$ belonging to $Acc^\eta
(T)\backslash \{(T,0,\ldots,0)\}$ and to the
$\xCzero$-neighborhood of $\gamma$ of Lemma \ref{formenormaleexcept},
we have~: $x_n>0$.

So Lemma \ref{lemma2} describes the $n-2$ first coordinates and
the last coordinate.
Anyway \cite{BK} do not control the coordinate $E_{n-1}^T(u)$.
Indeed their method consists in replacing the control $u\in \xLtwo$
by the control $x_{n-1}\in \xLtwo$, which is
called \it{Goh transformation} and consists in fact in immersing
$\xLtwo$ into the space $H^{-1}$, see \cite{AS1}. Such trajectories
are called \it{generalized trajectories}, and this corresponds to
add to the set of controls $\xLtwo([0,T])$ Dirac measures at $0$
and $T$.

Here, in order to take into consideration the constraint $|u|\leq
\eta$ and to control the coordinate $E_{n-1}(u)$,
the method used in \cite{BK} has to be \it{adapted}.
The aim is to describe the boundary of the intersection of
$Acc^\eta(T)$ with the
plane $(x_1,x_n)$, close to the point $(T,0)$ (corresponding to
$\gamma(T)$). Let $x$ be a real close to $T$.
We have
to solve equations~:
\begin{equation} \label{equas1}
E_1^T(u)=x,\ E_2^T(u)=\cdots=E_{n-1}^T(u)=0
\end{equation}
in the domain $|u|\leq\eta$, and then among such solutions we have
to minimize the last coordinate $E_n^T(u)$ (in order to describe
the boundary of the accessibility set).
We proceed in the following way. Set
$x_1(t)=t+\xi(t)$~; we shall first determine $\xi$ such that the
control $u=\xi^{(n-1)}$ satisfies (\ref{equas1}). From Lemma
\ref{lemma1}, $\xi$ can be expanded in~:
\begin{equation} \label{expansionxi}
\xi=\sum_{i=1}^{n-2}\alpha_iJ_i^T
+\sum_{i=1}^\infty v_ie_i^T
\end{equation}
From Lemma (\ref{lemma2}), the first $n-2$ coordinates are,
if $\ds{u\sim
\xi^{(n-1)}=\sum_{i=1}^{n-2}\alpha_i{J_i^T}^{(n-1)}+v}$~:
$$\ds{\begin{array}{rcl}
E_1^T(u)&=&T+\alpha_1+\rm{o}(\alpha_1,\ldots,\alpha_{n-2},
||v||_{\xLinfty}) \\
E_2^T(u)&=&\alpha_2+\rm{o}(\alpha_1,\ldots,\alpha_{n-2},
||v||_{\xLinfty}) \\
&\vdots& \\
E_{n-2}^T(u)&=&\alpha_{n-2}+\rm{o}(\alpha_1,\ldots,\alpha_{n-2},
||v||_{\xLinfty}) \\
\end{array}}$$
Let us solve equations (\ref{equas1}). We get from the
\it{Implicit Function Theorem}~:
$$\begin{array}{rcl}
\alpha_1&=&x-T+\rm{o}(x-T,||v||_{\xLinfty})\\
\alpha_2&=&\rm{o}(x-T,||v||_{\xLinfty})\\
&\vdots& \\
\alpha_{n-2}&=&\rm{o}(x-T,||v||_{\xLinfty})\\
\end{array}$$
as $x\rightarrow T$ and $v\overset{\xLinfty}{\longrightarrow}0$.
Further we will modify $u$ in order to get also~:
$$x_{n-1}(0)=x_{n-1}(T)=0\quad\rm{and}\quad |u|\leq\eta$$
Let us first study the last coordinate. We get from Lemma
(\ref{lemma2})~:
$$E_n^T(u)=A_T(x-T)^2+\sum_{i=1}^\infty \lambda_i^Tv_i^2
+\rm{o}\l( \l( x-T+||v||_{\xLinfty} \r)^2 \r)$$
as $x\rightarrow T$ and $v\overset{\xLinfty}{\longrightarrow}0$.
In this expression we can see that the minimum is reached at
$v=0$ \it{up to a $\rm{o}((x-T)^2)$}, for if $T<t_c$ then all
$\lambda_i^T$'s are positive (Lemma \ref{lemma0}).

\begin{rmrk}
If we do not neglect any term then the minimum is not
necessarily reached at $v=0$ because of terms of order $3$
in $(x-T)^2v_i$.
\end{rmrk}

\no Therefore, \it{at order $2$}, the minimum is reached at
$\ds{u\sim\sum_{i=1}^{n-1}\alpha_i{J_i^T}^{(n-1)}}$ and equals~:
$$E_n^T(u)=A_T(x-T)^2+\rm{o}((x-T)^2)\qquad \rm{ as }
x\rightarrow T$$

We shall now prove that it is possible to modify $u$, taking
into account the constraint $|u|\leq\eta$,
in order to obtain $x_{n-1}(0)
=x_{n-1}(T)=0$, without changing the previous results on the
other coordinates.

\no First of all, if $x$ is close enough
to $T$, then the $\alpha_i$'s are
small and thus the constraint $|u|\leq\eta$ is satisfied.

\no Let us modify $u$ (and hence $x_{n-1}$) in the following way.
Set $u(t)=\eta$ if $x_{n-1}(0)>0$ or $-\eta$ if $x_{n-1}(0)<0$
on $[0,t_1]$ (same construction on $[t_2,T]$),
where $t_1$ and $t_2$ are the coinciding times, i.e. the times
at which $x_{n-1}$ coincides
with its initial graph, see Fig. \ref{figgraphx}.

\drawingscale{0.35mm}
\placedrawing[h]{fgraphx.lp}{}{figgraphx}

Precisely~:
\begin{itemize}
\item if $0\leq t\leq t_1$~: $|u(t)|=\eta$, $|x_{n-1}(t)|=\eta t$,
and the coinciding time
$t_1$ is such that $\ds{\eta t_1\sim|\sum_{i=1}^{n-2}\alpha_i
{J_i^T}^{(n-2)}(t_1)| }$.
\item if $t_1\leq t\leq t_2$~: $u(t)\sim\xi^{(n-1)}(t)$,
$\ds{x_{n-1}(t)=\sum_{i=1}^{n-2}\alpha_i
{J_i^T}^{(n-2)}(t)  }$.
\item if $t_2\leq t\leq T$~: $|u(t)|=\eta$, $|x_{n-1}(t)|
=\eta(T-t)$, where $t_2$ is such that
$\ds{\eta(T-t_2)\sim|\sum_{i=1}^{n-2}\alpha_i{J_i^T}^{(n-2)}(t_2)|}$.
\end{itemize}

We shall now check that if
$$x-T=\rm{o}(\eta)\quad \rm{as }x\rightarrow T$$
then all previous results are still valid.

\no Clearly~:
$\eta t_1=\rm{O}(x-T)$ and $\eta(t_2-T)=\rm{O}(x-T)$ as
$x\rightarrow T$. Moreover, at order $1$, the system is~:
$$\dot{x}_1=1+x_2,\ \dot{x}_2=x_3,\ \ldots,\ \dot{x}_{n-1}=u$$
Therefore~:
$$x(t_1)-t_1\sim\eta\f{t_1^{n-1}}{(n-1)!} \quad\rm{and}\quad
x_1(t_2)-t_2-x_1(T)+T\sim\eta\f{t_2^{n-1}}{(n-1)!}$$
Now if $n\geq 3$ and $x-T=\rm{o}(\eta)$ then these terms are
negligibly small in comparison to $x-T$.

\no As concerns the last coordinate, we obtain~:
$$x_n(t_1)=\rm{O}(\eta^2t_1^3)=\rm{o}((x-T)^2) \quad\rm{and}\quad
x_n(T)-x_n(t_2)=\rm{O}(\eta^2t_2^3)=\rm{o}((x-T)^2)$$
In these conditions, all our previous construction is still
valid. Hence in the plane $(x_1,x_n)$ the boundary of
$Acc^\eta(T)$ is a curve of class $\xCtwo$, independant of the
constraint, such that $x_n\sim A_T(x_1-T)^2$, which proves
Theorem \ref{bordaffine}. Moreover it results
from \cite{BK} and \cite{Sar} that the function $T\mapsto A_T$ is
continuous and decreasing on $[0,t_c[$, positive on $[0,t_{cc}[$
and negative on $]t_{cc},t_c[$.
\end{proof}


\subsection{Application to the sub-Riemannian case}
\label{sectionsrcase}
\subsubsection{Asymptotics of the sub-Riemannian sphere along
an abnormal direction}
Consider a smooth
sub-Riemannian structure $(M,\Delta,g)$ where $M$ is a
Riemannian $n$-dimensional manifold, $n\geq 3$, $\Delta$ is
a rank 2 distribution on $M$, and $g$ is a metric on $\Delta$. Let
$x_0\in M$~; our point of view is local and we can assume that
$M=\xR^n$ and $x_0=0$. Suppose there exists a smooth injective
abnormal trajectory
$\gamma$ passing through $0$. Up to changing coordinates and
reparametrizing we can assume that~:
\begin{itemize}
\item $\gamma(t)=(t,0,\ldots,0)$,
\item $\Delta=\rm{Span }\{X,Y\}$ where $X,Y$ are $g$-orthonormal,
\item $\gamma$ is the integral curve of $X$ passing through $0$.
\end{itemize}
Under these assumptions, the sub-Riemannian problem is equivalent
to the \it{time-optimal problem} for the system~:
\begin{equation} \label{systemesousR}
\dot{x}=vX(x)+uY(x),\ x(0)=0
\end{equation}
where the controls $v,u$ satisfy the \it{constraint}~:
\begin{equation}
v^2+u^2\leq 1
\end{equation}
The reference abnormal trajectory $\gamma$ corresponds to the
control~: $v=1,u=0$.

Let us now define a notion of \it{constrained accessibility
set}~:

\begin{dfntn}
Let $0<\alpha<1$. We denote by $Acc_{SR}^\alpha(T)$ the
accessibility set at time $T$ for the sub-Riemannian
system (\ref{systemesousR})
with the additional constraint on the control~:
$$v^2+u^2\leq 1\ ,\   1-\alpha\leq v\leq 1\ ,\   |u|\leq
\alpha$$
(see Fig. \ref{figcontSR})
\end{dfntn}

\drawingscale{0.3mm}
\placedrawing[h]{fcontsr.lp}{}{figcontSR}

Note that controls steering $0$ to points of $Acc_{SR}^\alpha(T)$
are in a \it{$\alpha$-neighborhood} in $\xLinfty$ metric
of the abnormal reference control $v=1,u=0$.

\begin{dfntn}
We call \it{affine system associated to the sub-Riemannian
system} (\ref{systemesousR}) the following system~:
\begin{equation} \label{affi}
\dot{x}=X(x)+wY(w)
\end{equation}
where the control $w$ satisfies a constraint of the form~:
$|w|\leq \eta$.
\end{dfntn}

Let $Acc^\eta_A(T)$ denote the accessibility set at time $T$
for this affine system with the constraint~: $|w|\leq \eta$. The
reference trajectory $\gamma$ corresponds to $w=0$, and is also
\it{abnormal} for this affine system.

The following lemma gives a precise comparison of constrained
accessibility sets of systems (\ref{systemesousR}) and
(\ref{affi})~:

\begin{lmm} \label{lemmecomparaison}
\begin{enumerate}
\item $\ds{\forall \alpha\in ]0,1[\quad Acc_{SR}^\alpha(T)
\subset \bigcup_{(1-\alpha)T\leq s\leq T}
Acc_A^{\f{\alpha}{1-\alpha}}(s)}$ 
\item $\ds{\bigcup_{T_0\leq s\leq \f{T}{\sqrt{1+\eta^2}}}
Acc_A^\eta(s) \subset Acc_{SR}^\alpha(T) }\quad
\rm{ where } \alpha=\max \l(1-\f{T_0}{T},\f{\eta}
{\sqrt{1+\eta^2}}\r)$.
\end{enumerate}
\end{lmm}

\begin{proof}
Let us prove the first inclusion. If $x_1\in Acc_{SR}^\alpha(T)$
then there exists a control $(v,u)$ such that $v^2+u^2\leq 1$,
$1-\alpha\leq v\leq 1$, $|u|\leq \alpha$, and such that the
corresponding trajectory satisfies~:
$$\dot{x}=vf_0+uf_1,\ x(0)=0,\ x(T)=x_1$$
As $\alpha<1$, $v$ does not vanish, the following
reparametrizing holds~: $\f{ds}{dt}=v$. Set $y(s)=x(t)$,
$w(s)=\f{u(t)}{v(t)}$, and $S=\int_0^Tv$. Then~: $x_1=x(T)=y(S)$,
and $S$ is such that~:
$$(1-\alpha)T\leq S\leq T$$
Moreover~:
$$\f{dy}{ds}=f_0+\f{u}{v}f_1=f_0+wf_1$$
where $|w|\leq\f{\alpha}{1-\alpha}$. Therefore~:
$$x_1\in \bigcup_{(1-\alpha)T\leq s\leq T}
Acc_A^\f{\alpha}{1-\alpha}(s)$$
which proves the first part of the lemma.

Let us now check the second inclusion. Let $S\in
[T_0,\f{T}{\sqrt{1+\eta^2}}]$ and $x_1\in Acc_A^\eta(S)$. There
exists a control $w$ such that $|w|\leq \eta$ and
the corresponding trajectory satisfies~:
$$\f{dy}{ds}=f_0+wf_1,\ y(0)=0,\ y(S)=x_1$$
Let $\varepsilon >0$ such that $T=\f{S}{1-\varepsilon }$. Let us
make the reparametrizing~: $\f{ds}{dt}=1-\varepsilon $, and set~:
$x(t)=y(s)$, $v(t)=1-\varepsilon $, $u(t)=(1-\varepsilon )w(s)$,
where $t\in[0,T]$. Then~:
$$\dot{x}=vf_0+uf_1,\ x(0)=0,\ x(T)=y(S)=x_1$$
Let us now check the constraint on the control $(v,u)$. By
definition~: $T_0\leq S\leq \f{T}{\sqrt{1+\eta^2}}$ and
$S=(1-\varepsilon )T$. Hence~:
$$1-\varepsilon \leq \inv{\sqrt{1+\eta^2}}\quad\rm{ and }\quad
1-\varepsilon \geq \f{T_0}{T}$$
and thus~:
$$1-\l(1-\f{T_0}{T}\r)\leq v\leq 1\quad\rm{ and }\quad
|u|\leq\f{\eta}{\sqrt{1+\eta^2}}$$
Moreover~:
$$v^2+u^2\leq (1-\varepsilon )^2(1+\eta^2)\leq 1$$
Therefore we can conclude that~:
$$x_1\in Acc_{SR}^\alpha(T)\quad\rm{ where }\quad
\alpha=\max\l(1-\f{T_0}{T},\f{\eta}{\sqrt{1+\eta^2}}\r)$$
\end{proof}

\no Using the previous lemma and Theorem
\ref{bordaffine} we can prove the following~:

\begin{thrm} \label{bordSR}
Suppose assumptions ($H_0-H_3$) are fulfilled along the reference
abnormal trajectory $\gamma$ for the system $(X,Y)$. Let
$t_{cc}$ and $t_c$ denote the first
conjugate times of $\gamma$ for the
associated affine system. Let $\alpha\in ]0,1[$. Then~:
\begin{enumerate}
\item There exist coordinates $(x_1,\ldots,x_n)$ locally along
$\gamma $ such that in these coordinates~:
$\gamma(t)=(t,0,\ldots,0)$, and the first Pontryagin's cone along
$\gamma $ is~: $K(t)=\rm{Vect
}\{\f{\partial}{\partial x_1},\ldots,\f{\partial}{\partial
x_{n-1}} \}_{|\gamma}$.
\item If $T$ is small enough then for any point
$(x_1,\ldots,x_n)$ of $Acc^\alpha_{SR}(T)$ close to $\gamma(T)$ we
have $x_n\geq 0$ (see Fig. \ref{figparab}).
\item If $T<t_{cc}$, then in the plane $(x_1,x_n)$, close to the
point $(T,0)$, the boundary of $Acc_{SR}^\alpha(T)$ does not
depend on $\alpha$, is a curve of class $\xCtwo$ outside $(T,0)$,
tangent to the abnormal direction, whose first term is~:
\begin{itemize}
\item if $x_1\leq T$ then $x_n=0$.
\item if $x_1\geq T$ then $x_n=A_T(x_1-T)^2+\rm{o}((x_1-T)^2)$.
\end{itemize}
The function $T\mapsto A_T$ is the same as in Theorem
\ref{bordaffine}.
\item If $T>t_{cc}$ then $Acc^\alpha_{SR}(T)$ is open near
$\gamma(T)$.
\end{enumerate}
\end{thrm}

Figure \ref{figbordSR} represents the evolution of
$Acc^\alpha_{SR}(T)$ in function of $T$ in the plane
$(x_1,x_n)$. It is open in a
neighborhood of $\gamma(T)$ if $T>t_{cc}$, contrarily to the
affine case where it becomes open only beyond $t_c$.

\drawingscale{0.5mm}
\placedrawing[h]{fbordsr.lp}{}{figbordSR}

\begin{rmrk}
To compare the system (\ref{systemesousR}) with its associated
affine system (\ref{affi}) we need the following reparametrizing
(see Proof of Lemma \ref{lemmecomparaison})~:
$$\f{ds}{dt}=v$$
which only holds if $v$ does not vanish. This condition is
satisfied when the control $(v,u)$ is in a
$\alpha$-neighborhood in $\xLinfty$ metric
of the abnormal reference control
$(1,0)$, for in this case $v$ is close to $1$ in $\xLinfty$.
Hence using this method it is only possible to describe a
constrained accessibility set, i.e. in a $\alpha$-neighborhood in
$\xLinfty$ metric of the reference abnormal control.
\end{rmrk}

\begin{proof}[Proof of Theorem \ref{bordSR}]
The aim is to compare precisely systems (\ref{systemesousR}) and
(\ref{affi}) using Lemma \ref{lemmecomparaison} and to apply
Theorem \ref{bordaffine}. In order to do this we first have to
\it{normalize} the affine system (\ref{affi}) using
Lemma \ref{formenormaleexcept}. We denote by $Acc_A^\beta(T)$ the
accessibility set at time $T$ with constraint $\beta$ of the
affine system (\ref{affi}).

The system (\ref{systemAN}), $\dot{x}=f_0+uf_1$,
is called \it{normalized affine
system} and will be refered as (AN). Let $Acc_{AN}^\eta(T)$
denote the accessibility set at time $T$ for this system with the
constraint~: $|u|\leq\eta$. Due to the particular forms of the
feedback $(\varphi,\alpha,\beta)$ and of the system, we have~:

\begin{lmm} \label{lem56}
$\ds{\varphi ^{-1}(Acc_{AN}^\eta(T))\ =\ Acc_A^\beta(T)}$ where
$\beta=\rm{O}(\eta)$ as $\eta\rightarrow 0$.
\end{lmm}

\no We know from Lemma \ref{lemmecomparaison} that~:
$$Acc_{SR}^\alpha(T)\subset\bigcup_{s\leq T}Acc_{A}^
\f{\alpha}{1-\alpha}(s)$$
Hence in the normalized coordinates we get~:
\begin{equation} \label{etoile}
\varphi (Acc_{SR}^\alpha(T))\subset \bigcup_{s\leq
T}Acc_{AN}^\beta(s)
\end{equation}
We shall use Theorem \ref{bordaffine}, which describes the
boundary of
$Acc_{AN}^\beta(s)$ in the plane $(x_1,x_n)$, to study
the boundary of $\ds{\bigcup_{s\leq
T}Acc_{AN}^\beta(s)}$. Using the fact that
t The function $t\mapsto A_T$ is continuous and decreasing on
$[0,t_{cc}[$, we can assert that the boundary of
$\ds{\bigcup_{s\leq T}Acc_{AN}^\beta(s)}$ is given, in the plane
$(x_1,x_n)$, close to the point $(T,0)$, by the following curve,
see Fig. \ref{figbdproof}~:
\begin{itemize}
\item if $x_1\leq T$ then $x_n=0$.
\item if $x_1\geq T$ then $x_n=A_T(x_1-T)^2+\rm{o}((x_1-T)^2)$
\end{itemize}

\drawingscale{0.4mm}
\placedrawing[h]{fbdproof.lp}{}{figbdproof}

Let $x_n=f(x_1)$ denote the curve parametrizing the boundary of
$\varphi (Acc_{SR}^\alpha(T))$ in the plane $(x_1,x_n)$. From
inclusion (\ref{etoile}) we get~:
\begin{itemize}
\item if $x_1\leq T$ then $f(x_1)\geq 0$,
\item if $x_1\geq T$ then $f(x_1)\geq A_T(x_1-T)^2+\rm{o}
((x_1-T)^2)$.
\end{itemize}

Let us now prove the converse inequality. To this aim we shall
use varying constraints depending on $x_1-T$. We proceed in the
following way.
From Lemma \ref{lem56}~:
$$\varphi ^{-1}\Big(Acc_{AN}^\eta\Big(\f{T}{\sqrt{1+\eta^2}}
\Big)  \Big)
\subset Acc_A^\beta\Big(\f{T}{\sqrt{1+\eta^2}}\Big)$$
where $\beta=\rm{O}(\eta)$. Without loss of generality we can
assume $\beta\leq\eta$, and thus~:
$$Acc_A^\beta\Big(\f{T}{\sqrt{1+\eta^2}}\Big) \subset
Acc_A^\eta\Big(\f{T}{\sqrt{1+\eta^2}}\Big)$$
Now from Lemma \ref{lemmecomparaison}, we get that for any
$\eta>0$ small enough~:
$$Acc_A^\eta\Big(\f{T}{\sqrt{1+\eta^2}}\Big)\subset
Acc_{SR}^\alpha(T)$$
And thus in the normalized coordinates~:
\begin{equation} \label{etoile2}
Acc_{AN}^\eta\Big(\f{T}{\sqrt{1+\eta^2}}\Big)\subset
\varphi (Acc_{SR}^\alpha(T))
\end{equation}
Let $f_{AN}(x_1)$ be the function parametrizing the boundary of
$Acc_{AN}^\eta\Big(\f{T}{\sqrt{1+\eta^2}}\Big)$ in the plane
$(x_1,x_n)$. We know from Theorem \ref{bordaffine} that~:
$$f_{AN}(x_1)=A_{\f{T}{\sqrt{1+\eta^2}}}\Big(x_1-
\f{T}{\sqrt{1+\eta^2}}\Big)^2
+\rm{o}\Big(\Big(\f{T}{\sqrt{1+\eta^2}}
\Big)^2 \Big)$$
provided $\ds{x_1-\f{T}{\sqrt{1+\eta^2}}=\rm{o}(\eta)}$. This
latter condition is fulfilled if $\eta=(x_1-T)^\f{3}{4}$, and in
this case we have moreover~: $\eta^2=\rm{o}(x_1-T)$.
On the other part, from the continuity of $t\mapsto A_t$~:
$$A_{\f{T}{\sqrt{1+\eta^2}}}=A_T+\rm{o}(1)\qquad \rm{as
}x_1\rightarrow T$$
We obtain~: $f_{AN}(x_1)=A_T(x_1-T)^2+\rm{o}((x_1-T)^2)$.
Finally, from inclusion (\ref{etoile2}) we conclude~:
$$f(x_1)\leq A_T(x_1-T)^2+\rm{o}((x_1-T)^2)$$
which ends the proof.
\end{proof}


\subsubsection{Splitting-up of the sphere near an abnormal
direction} \label{secteursSR}
Let $T>0$ small enough so that properties 2 and 3
of Theorem \ref{bordSR}
are satisfied. In particular the reference abnormal trajectory
$\gamma$ is minimizing, see Section \ref{sectionoptimalitySR}.
Then $A=\gamma(T)$ belongs to the
sub-Riemannian sphere $S(0,T)$ with radius $T$.
If controls steering $0$ to points of the boundary of
$Acc^\alpha_{SR}(T)$ in $x_n>0$
(that are \it{$\xLinfty$-optimal}) are
actually \it{globally optimal}, then this boundary is included in
the sphere $S(0,T)$. In this case the sphere splits into \it{two
sectors} near $\gamma(T)$, bordered by the first Pontryagin's cone
$x_n=0$~:
\begin{itemize}
\item sector $x_n>0$ corresponding to the previous description,
\item sector $x_n<0$.
\end{itemize}
According to the previous results, final points at time $T$
associated to controls which are $\xLinfty$-close to the
reference abnormal
control are in the first sector~: $x_n>0$. Obviously
due to controllability of the system the sector $x_n<0$ is
accessible. In fact a basic calculus shows~:

\begin{lmm} \label{lem59}
For any neighborhood $V$ of the point $A$ in $\xR^n$ we have~:
$$S(0,T)\cap V\cap (x_n<0) \neq \emptyset$$
\end{lmm}

These points in $(x_n<0)$ are reached by controls which are close
to the reference control \it{in $\xLtwo$ metric but not in
$\xLinfty$ metric}. More precisely~:

\begin{lmm} \label{lem58}
Let $M_n=E_T(u_n)\in S(0,T)$ whose last coordinate $x_n$ is
strictly negative. Let $u$ denote the abnormal reference control.
We suppose that $M_n$ converges to $A= E(u)$ in $\xR^n$.
Then $u_n$ converges to
$u$ in $\xLtwo([0,T])$ but not in $\xLinfty([0,T])$.
\end{lmm}

Hence near the abnormal direction the sphere is splits into two
sectors~:
\begin{itemize}
\item the \it{$\xLinfty$-sector}~: $(x_n>0)\cap S(0,T)$
(described by Theorem \ref{bordSR}), made of end-points of
minimizing trajectories associated to controls $\xLinfty$-close
to the abnormal reference control,
\item the \it{$\xLtwo$-sector}~: $(x_n<0)\cap S(0,T)$, made of
points reached by minimizing controls $\xLtwo$-close, but
not $\xLinfty$-close to the abnormal reference control.
\end{itemize}
The contact of
the first sector is known, but not the second one a priori.
Anyway according to the \it{Tangency Theorem} (see \cite{T}),
under some nice stratification assumptions, this $\xLtwo$-sector
\it{ramifies tangently} to the Pontryagin cone $x_n=0$, see Fig.
\ref{figsect}.

\drawingscale{0.5mm}
\placedrawing[h]{fsect.lp}{}{figsect}

\begin{rmrk}
In particular, minimizing trajectories $\gamma_\infty$
joining $0$ to
points of the sphere in the $\xLinfty$-sector are
\it{$\xCone$-close} to the reference abnormal trajectory
$\gamma$. Minimizing trajectories $\gamma_2$
joining $0$ to points of the sphere in the
$\xLtwo$-sector are \it{$\xCzero$-close, but not $\xCone$-close} to
$\gamma$, see Fig. \ref{figsect}.
\end{rmrk}

\paragraph{{\bf{Typical example~: the Martinet case.}}}
Consider the two following vector fields in $\xR^3$~:
$$X=\f{\partial}{\partial x}+\f{y^2}{2}\f{\partial}{\partial z}
\ ,\   Y=\f{\partial}{\partial y}$$
and endow the distribution spanned by these vector fields with an
analytic metric $g$ of the type~:
$$g=adx^2+cdy^2$$
where $a=(1+\alpha y)^2$ and $c=(1+\beta x+\gamma y)^2$. The
abnormal reference control for the sub-Riemannian system
$\dot{x}=vX(x)+uY(x)$ with constraint $v^2+u^2\leq 1$ is
$v=1,u=0$, and corresponds to the trajectory $\gamma$~:
$x(t)=t, y(t)=z(t)=0$. We have, see \cite{BT}~:

\begin{lmm} \label{bordSRMartinet}
Assumptions ($H_0-H_3$) are fulfilled along $\gamma$
if and only if
$\alpha \neq 0$. In this case branches $1$ and $2$ (see Fig.
\ref{figsect} with $x_1=x, x_n=z$) have the following contacts
with the abnormal direction~:
\begin{itemize}
\item branch $1$~: $x\geq T,\ 
z=\inv{2T\alpha^2}(x-T)^2+\rm{o}((x-T)^2)$
\item branch $2$~: $x\leq T,\ z\sim \inv{6}(1+\rm{O}(T))(x-T)^3$
\end{itemize}
\end{lmm}

\begin{rmrk}
The coefficient $A_T$ of the first branch can be computed
directly or using formula
(\ref{calculA_T}) (see Remark \ref{remcalculA_T}).

As we are in dimension 3, results of Theorem \ref{bordSR} are in
fact available on $\xR^+$, see Remark \ref{remtcc}. The
$\xLtwo$-sector is $z<0$ and the $\xLinfty$-sector is $z>0$.
\end{rmrk}

\begin{proof}[Proof of Lemma \ref{lem59}]
Consider the sub-Riemannian system (\ref{systemesousR}) in the
coordinates $(x_1,\ldots,x_n)$ of Theorem \ref{bordSR}~:
$$\begin{array}{rcl}
\dot{x}_1&=&v(1+x_2+x_nf_1(x_1)+Z_1)\\
\dot{x}_2&=&v(x_3+x_nf_2(x_1)+Z_2)\\
&\vdots &\\
\dot{x}_{n-1}&=&u+x_nf_{n-1}(x_1)+Z_{n-1}\\
\ds{\dot{x}_n}&=&\ds{v\sum_{i=2}^na_{ij}(x_1)x_ix_j+Z_n}
\end{array}$$
where the 1-jet (resp. the 2-jet) of $Z_1,\ldots,Z_{n-1}$ (resp.
$Z_n$) along $u=0$ is equal to $0$.

Our goal is to construct a control close to the reference
abnormal control in $\xLtwo$ metric, whose associated
final point at time $T$
is close to $(T,0,\ldots,0)$ and is such that
$x_n<0$. Let $\varepsilon >0$. Consider the following control
$v$ (see Fig. \ref{figv})~:
\begin{itemize}
\item if $0\leq t\leq \f{T-\varepsilon }{2}$ then
$v(t)=1$,
\item if $\f{T-\varepsilon }{2}
\leq t\leq \f{T+\varepsilon }{2}$ then
$v(t)=-1$,
\item if $\f{T+\varepsilon }{2}\leq t\leq T$ then
$v(t)=1$.
\end{itemize}

\drawingscale{0.3mm}
\placedrawing[h]{fv.lp}{}{figv}

Set $u=0$. It is clear that $(v,u)$ is abnormal (but not
minimizing). Consider the following perturbation (see Fig.
\ref{figdx})~:
\begin{itemize}
\item if $0\leq t\leq \f{T-\varepsilon }{2}$ then $\delta
v(t)=\delta u(t)=0$,
\item if $\f{T-\varepsilon }{2}< t\leq \f{T}{2}$ then $\delta
v(t)=1-\sqrt{1-\varepsilon ^2},\ \delta u(t)=\varepsilon $,
\item if $\f{T}{2}<t\leq\f{T+\varepsilon }{2}$ then $\delta
v(t)=1-\sqrt{1-\varepsilon ^2},\ \delta u(t)=-\varepsilon $,
\item if $\f{T+\varepsilon }{2}<t\leq T$ then $\delta v(t)=\delta
u(t)=0$.
\end{itemize}
It is clear that $(v+\delta v)^2+(\delta u)^2=1$. Moreover~:
$||\delta v||_{\xLone}=\rm{O}(\varepsilon ^3)$ and
$||\delta u||_{\xLone}=\rm{O}(\varepsilon ^2)$.

\drawingscale{0.3mm}
\placedrawing[h]{fdx.lp}{}{figdx}

The end-point mapping $E_T$ being $\xCinfty$ in $\xLone$ topology
(see for instance \cite{LK}), we have~:
\begin{equation} \label{ref1}
\forall i\in\{1,\ldots,n-1\}\quad
E_T^i(v+\delta v,\delta u)=E_T^i(v,0)+dE_T^i(v,0).(\delta
v,\delta u)+\rm{O}(\varepsilon ^4)
\end{equation}
and
\begin{equation} \label{ref2}
E_T^n(v+\delta v,\delta u)=E_T^n(v,0)+dE_T^n(v,0).(\delta
v,\delta u)+\inv{2}d^2E_T^n(v,0).(\delta v,\delta u)^{[2]}
+\rm{O}(\varepsilon ^6)
\end{equation}
Moreover we have $E_T^i(v,0)=\delta _i^1(T-\varepsilon )$ for
$i=1,\ldots,n$. On the other part the linearized system along
$(v,u)$ is~:
$$\dot{y}_1=\delta v+vy_2,\
\dot{y}_2=vy_3,\
\ldots,\
\dot{y}_{n-2}=vy_{n-1},\
\dot{y}_{n-1}=\delta u,\
\dot{y}_n=0$$
Let us calculate $y_{n-1}(t)$~:
\begin{itemize}
\item if $0\leq t\leq \f{T-\varepsilon }{2}$ then $y_{n-1}(t)=0$,
\item if $\f{T-\varepsilon }{2}< t\leq \f{T}{2}$ then
$y_{n-1}(t)=\varepsilon (t-\f{T-\varepsilon }{2})$,
\item if $\f{T}{2}<t\leq\f{T+\varepsilon }{2}$ then
$y_{n-1}(t)=-\varepsilon (t-\f{T+\varepsilon}{2})$,
\item if $\f{T+\varepsilon }{2}<t\leq T$ then $y_{n-1}(t)=0$.
\end{itemize}
In particular~: $||y_{n-1}||_{\xLinfty}=\rm{O}(\varepsilon ^2)$,
and $y_{n-1}$ is equal to $0$ outside an interval of length
$\varepsilon $. We get easily~:
\begin{equation} \label{ref3}
||y_{n-2}||_{\xLinfty}=\cdots=||y_1||_{\xLinfty}=\rm{O}(\varepsilon
 ^3)
\end{equation}
Hence using (\ref{ref1})~:
$$\begin{array}{rcl}
E^T_1(v+\delta v,\delta u)&=&T-\varepsilon +\rm{O}(\varepsilon ^3)
\\
E^T_2(v+\delta v,\delta u)&=&\rm{O}(\varepsilon ^3)
\\
& \vdots &\\
E_{n-2}^T(v+\delta v,\delta u)&=&\rm{O}(\varepsilon^3) \\
E^T_{n-1}(v+\delta v,\delta u)&=&\rm{O}(\varepsilon ^2)
\end{array}$$
Let us now compute the last coordinate. We have to calculate
$z_n$, where~:
$$\dot{z}_n=v\sum_{i,j=2}^{n-1}a_{ij}(t)y_iy_j$$
It represents the intrinsic second-order derivative. From
(\ref{ref3}) we get~:
$$\dot{z}_n(t)=v(t)a_{n-1,n-1}(t)y^2_{n-1}(t)+y_{n-1}(t)
\rm{O}(\varepsilon ^3)+\rm{O}(\varepsilon ^6)$$
As $y_{n-1}$ is equal to $0$ outside $[\f{T-\varepsilon
}{2},\f{T+\varepsilon }{2}]$, we have~:
$$\int_0^Ty_{n-1}(t)\rm{O}(\varepsilon ^3)dt\ =\
\rm{O}(\varepsilon ^6)$$
Moreover the coefficient $a_{n-1,n-1}$ is continuous and does not
vanish on $[0,T]$, hence there exists $\alpha>0$ such that
$a_{n-1,n-1}(t)\geq\alpha$ on $[0,T]$. We get~:
$$z_n(T)\leq -\alpha\int_{\f{T-\varepsilon
}{2}}^{\f{T+\varepsilon }{2}}y_{n-1}^2(t)dt+\rm{O}(\varepsilon
^6) \leq -\f{\alpha\varepsilon ^5}{12}+\rm{O}(\varepsilon ^6)$$
Therefore from (\ref{ref2})~:
$$E^T_n(v+\delta v,\delta u)\leq -\f{\alpha\varepsilon
^5}{12}+\rm{O}(\varepsilon ^6)$$
Hence for any neighborhood $V$ of $(T,0,\ldots,0)$, $Acc(T)\cap
V$ contains points such that $x_n<0$, and hence the same goes for
$S(0,T)\cap V$
since the abnormal reference trajectory is minimizing. By
construction, controls steering to these points are $\xLtwo$-close
but not $\xLinfty$-close to the reference abnormal control.
\end{proof}

\begin{proof}[Proof of Lemma \ref{lem58}]
As $M_n$ belongs to the sphere $S(0,T)$, we have~:
$||u_n||_{\xLtwo}=||u||_{\xLtwo}=T$, thus the sequence $(u_n)_{n\in\xN}$
is bounded in $\xLtwo$. Hence up to a subsequence we can assume that
$u_n$ converges weakly to $v\in \xLtwo$ (denoted by
$u_n\hookrightarrow v$). From the continuity of the
end-point mapping $E_T$ in the weak topology on
$\xLtwo$ (see \cite{T}), we can assert that $M_n=E_T(u_n)$ converges
to $E_T(v)$, and thus $E_T(v)=E_T(u)$. The assumptions on the
reference abnormal trajectory imply that $v=u$. Hence
$u_n\hookrightarrow u$, and on the other part~:
$||u_n||_{\xLtwo}=||u||_{\xLtwo}$, therefore $u_n$ converges (strongly)
towards $u$ in $\xLtwo$.

Moreover results stated by Theorem \ref{bordSR} imply that $u_n$
does not converge towards $u$ in $\xLinfty$ (because it would
imply that $x_n\geq 0$).
\end{proof}


\section{Application : optimality of abnormal trajectories}
\label{sectionoptimalité}
In this Section we apply our previous theory on accessibility
sets to studying optimality of abnormal trajectories~; this leads
us to find again some well-known results. Indeed
in the notations of Theorem \ref{bordaffine}, this Theorem
implies in particular that if $T<t_{cc}$ then $\gamma$ is
\it{isolated in $\xCzero$-topology}
in the space of all trajectories which
connect given end-points, and thus is \it{optimal for any cost}
in this topology. This well-known
property, called \it{rigidity}, was intensively studied. The
main results concerning this analysis in a generic context were
given first in \cite{Sar} and \cite{BK} for
single-input affine systems, then in \cite{AS,BHsu,
LS,Zhong}, for
sub-Riemannian systems, and in \cite{AS2} in general. Moreover
these authors developed a Morse theory in order to characterize
conjugate points, that is, points beyond which the abnormal
trajectory is no more optimal.

Hence results given in this Section are not really new.
However they are slightly different from the results cited above.
Indeed on the one part in \cite{BK} were obtained necessary and
sufficient conditions for $\xCzero$-time-optimality of abnormal
trajectories of single-input affine systems without any
constraint on the control. Here we improve their statement by
adding a \it{constraint on
the control} and studying the problem of minimizing \it{any
cost}. On the other part, in \cite{AS2} was made a general theory
(i.e. for nonlinear systems) on optimality
of abnormal trajectories in $\xLinfty$ topology on the controls.
Results given here are valid
in the $\xCzero$-topology on the trajectories (but only for
single-input affine systems).
Moreover we study the equivalence between the time-optimality
problem and the problem of minimizing any cost, the final time
being fixed or not.
Finally Theorem \ref{thrmoptimaliteSR}, which concerns
optimality of abnormal trajectories for sub-Riemannian systems of
rank 2, makes a
link between the works of \cite{BK} and \cite{AS,AS1}.

\subsection{Optimality of abnormal trajectories for single-input
affine systems} \label{sectionoptimalitySR}
Consider the single-input affine system (\ref{af}) with
constraint (\ref{contrainte}), and suppose assumptions
$(H_0-H_4)$ are fulfilled along a reference abnormal trajectory
$\gamma$.
We first study the time-optimal problem, and then the problem of
minimizing some cost.

\subsubsection{Time optimality} \label{sectiontimeoptimality}
\begin{dfntn}
\begin{itemize}
\item
The trajectory $\gamma$ is said \it{$\xCzero$-time-minimal} on
$[0,T]$ if there exists a $\xCzero$-neighborhood of $\gamma$ such
that $T$ is
the minimal time to steer $\gamma(0)$ to $\gamma(T)$ among the
solutions of the system (\ref{af}) with the constraint
(\ref{contrainte}) that are entirely contained in this
neighborhood.
\item Recall that $\gamma$ is associated to the control $u=0$. Let
$\delta>0$.
The trajectory $\gamma$ is said \it{$\xLinfty$-time-minimal} on
$[0,T]$ if there exists a neighborhood of $0$ in
$\xLinfty([0,T+\delta])$ such that $T$ is the minimal time to
steer $\gamma(0)$ to $\gamma(T)$ among trajectories associated to
controls of this neighborhood.
\end{itemize}
Obviously if $\gamma$ is $\xCzero$-time-minimal then it is
$\xLinfty$-time-minimal.
\end{dfntn}

\begin{thrm} \label{tempsminhyp}
Under assumptions of Theorem \ref{bordaffine},
the trajectory $\gamma$ is $\xCzero$-time-minimal if and only if
$T<t_{cc}$. Moreover $\gamma$
is not $\xLinfty$-time-minimal if $T>t_{cc}$.
\end{thrm}

The proof is clear by inspecting Fig. \ref{figaccexc} and Proof
of Theorem \ref{bordaffine}.

\begin{rmrk} \label{remtcc2}
If $n=3$, we have $t_{cc}=+\infty$ provided assumptions
($H_0-H_3$) are fulfilled on $\xR^+$. Hence in this case $\gamma$
is $\xCzero$-time-minimal on $\xR^+$.
\end{rmrk}


\subsubsection{Optimality for some cost}
Let us now consider the problem of minimizing some cost
$C(T,u)$, also denoted by $C_T(u)$, where $C$ is a smooth
function satisfying the following additional assumption along the
reference singular trajectory $\gamma$~:
\begin{itemize}
\item[($H_5$)] $\forall T\quad \rm{rank }(dE_T(0),dC_T(0))=n$
\end{itemize}
i.e. the singularity of the end-point mapping of the extended
system has codimension 1, and in particular the cost is independant
from the end-point mapping along $\gamma$.
We consider several optimization
problems~:
\begin{enumerate}
\item final time not fixed : the aim is to steer the system from
$x_0$ to $x_1$ in some time $T$ (not preassigned) and minimizing
the cost $C$.
\item final time fixed : let $T>0$ fixed~;
the aim is to steer the system
from $x_0$ to $x_1$ in time $T$ and minimizing the cost $C_T$.
\end{enumerate}

\paragraph{{\it{1. Final time not fixed}}}
\begin{dfntn}
\begin{itemize}
\item
The trajectory $\gamma$ is said to be \it{$\xCzero$-cost-minimal}
on $[0,T]$ if there exists a $\xCzero$-neighborhood of $\gamma$ such
that for any trajectory $q$ contained in this neighborhood, with
$q(0)=\gamma(0)$ and $q(t)=\gamma(T)$, we have~:
$C(t,v)\geq C(T,0)$,
where $v$ is the control associated to $q$.
\item Let $\delta>0$.
The trajectory $\gamma$ is said to be
\it{$\xLinfty$-cost-minimal} on $[0,T]$
if there exists a neighborhood of $0$ in $\xLinfty([0,T+\delta])$
such that, for any trajectory $q$ associated to a control $v$ of
this neighborhood, with $q(0)=\gamma(0)$ and
$q(t)=\gamma(T)$, we have~:
$C(t,v)\geq C(T,0)$.
\end{itemize}
Obviously the $\xCzero$-cost-minimality implies the
$\xLinfty$-cost-minimality.
\end{dfntn}

\no We have the following result (compare with \cite{AS2})~:

\begin{thrm} \label{coutmintfnf}
Under assumptions $(H_0-H_5)$,
the trajectory $\gamma$ is
$\xCzero$-cost-minimal if and only if it is
$\xCzero$-time-minimal. Actually,
$\gamma $ is $\xCzero$-cost-minimal if $T<t_{cc}$, and is not
$\xLinfty$-cost-minimal if $T>t_{cc}$.
\end{thrm}

Hence if the final time is not fixed then
\it{both problems of cost-minimization
and time-minimization are equivalent}.

\begin{proof}
If $T<t_{cc}$, then
the abnormal trajectory $\gamma$ is isolated in a
$\xCzero$-neighborhood, hence in particular is $\xCzero$-cost-minimal. If
$T>t_{cc}$, we know that $\gamma$ is no more time-minimal
since the point $(T,0,\ldots,0)$ belongs to $Acc(t)$ for some
$t<T$, see fig. \ref{fig32}.

\drawingscale{0.35mm}
\placedrawing[h]{f32.lp}{}{fig32}

To prove that $\gamma$ is no more $\xLinfty$-cost-minimal, we
have to show that \it{the point $(T,0,\ldots,0)$ belongs to
$M_r(t)$ for some $r<C_T(0)$}, where $C_T(0)$ is the cost of
$\gamma$ at time $T$, and $M_r(t)$ is the level set $r$ at time
$t$ of the value function associated to the cost $C$.

In the coordinates $(x_1,\ldots,x_n)$ of Lemma
\ref{formenormaleexcept}, the end-point mapping is
$E_t=(E_1^t,\ldots, E_n^t)$, where~:
$$\rm{rank }(dE^t_1(0),\ldots,dE^t_n(0))=n-1\quad\rm{ and }\quad
dE^t_n(0)=0$$
Instead of working in $\xLtwo([0,T])$, we choose as control space
the Sobolev space $\xHone([0,T])=\{f\in \xLtwo([0,T])\ /\ f'\in
\xLtwo([0,T])\}$. Endowed with the norm
$||f||_{\xHone}=\sqrt{||f||_{\xLtwo}^2 +||f'||_{\xLtwo}^2}$, this
is a Hilbert space. The reason to use it is the following~:
it can be immersed in a compact way in $\xCzero([0,T])$, and
in the proof we shall need continuous controls. 

The end-point mapping
is still differentiable in $\xHone$ since it is differentiable in
$\xLtwo$. From assumption ($H_4$), the end-point mapping is
independant from the cost along $\gamma$, hence at time $T$ there
exist $n$ independant vector fields $e_0^T,e_1^T,\ldots,
e_{n-1}^T$ such that~:
$$dC_T(0).e_0^T=1\quad\rm{ and }\quad
dE_i^T(0).e_i^T=1,\ i=1\ldots n-1$$
We can decompose $\xHone$ in the following way~:
$$\xHone([0,T])=\xR e_1^T\oplus \xR e_2^T\oplus\cdots\oplus\xR
e_{n-1}^T\oplus\xR e_0^T\oplus F$$
where $\xker dE_i^T(0)=\underset{j\neq i}{\bigoplus}\xR
e_j^T\oplus F$
and $\xker dC_T(0)=\underset{j\neq 0}{\bigoplus}\xR
e_j^T\oplus F$.

Let $Q_T(u)=d^2E_n^T(0).(u,u)$ (as in Lemma \ref{lemma2}).
The kernel of the first derivative is~: $\xker dE_T(0)=\xR
e_0^T\oplus F$, and hence the intrinsic second-order derivative
at time $T$ along $\gamma$ is the quadratic form ${Q_T}_{/\ \xR
e_0^T\oplus F}$. On the other part the space $\xR e_1^T\oplus\xR
 e_0^T\oplus F $ represents the domain of the latter quadratic
form where {the limit condition on the first coordinate has been
relaxed}, hence it is the domain of the \it{reduced operator
$D_1$}.
Now by definition of the first conjugate times $t_{cc}$
and $t_c$, we have~:
\begin{itemize}
\item if $T<t_{cc}$~: ${Q_T}_{/\ \xR e_1^T\oplus\xR e_0^T\oplus F}$
is positive definite,
\item if $t_{cc}<T<t_{c}$~: ${Q_T}_{/\ \xR e_0^T\oplus F}$
is positive definite, and ${Q_T}_{/\ \xR e_1^T\oplus\xR e_0^T\oplus
F}$ is indefinite,
\item if $T>t_c$~: ${Q_T}_{/\ \xR e_0^T\oplus F}$ is indefinite.
\end{itemize}
From now on we assume that $T>t_{cc}$. Then there exist
$v,w$ in $\xR e_0^T\oplus F$ such that $Q_T(e_1^T+v)<0$ and
$Q_T(w)>0$.

\begin{rmrk}
Note that we do not have necessarily $Q_T(e_1^T)<0$.
\end{rmrk}

Let $\delta>0$ small. As $v$ and $w$ are continuous, it is
possible to extend them on $[T,T+\delta]$ respectively by $v(T)$
and $w(T)$ (here is the role of $\xHone$). Consider now the
following control $u\in \xHone([0,T])$~:
$$u=a_1e_1^T+a_2e_2^T+\cdots+a_{n-1}e_{n-1}^T+a_1v+a_nw$$
where $a_1,\ldots,a_n$ are real numbers. In the same way we can
extend $u$ on $[0,T+\delta]$. We shall prove that we can choose
$a_1,\ldots,a_n$ small and $t$ close to $T$ such that the
previous control $u$ satisfies~:
$$|u|\leq\eta\ ,\   t\leq T+\delta\ ,\  
E_t(u)=\gamma(T)\ ,\   C_t(u)<C_T(0)$$
i.e. $\gamma(T)$ belongs to $M_r(t)$ for some $r<C_T(0)$. To this
aim we shall first use the \it{Implicit Function Theorem} on
$a_1,\ldots,a_n$ (to take into account the $n-1$ first
coordinates, i.e. the controllable part of the system),
then the \it{Mean Value Theorem} on the last
coordinate, and endly choose $t$ in order to make the cost lower
than $C_T(0)$. Note that in this method all reasonings are in
finite dimension.

Let us expand, for $t$ close to $T$ and $i=1,\ldots,n-1$~:
$$E_i^t(u)=E_i^t(0)+dE_i^t(0).u+\rm{o}(||u||_{\xHone([0,t])})$$
Moreover~:
\begin{itemize}
\item $E_1^t(0)=t$, and for $i=2,\ldots,n-1$~: $E_i^t(0)=0$,
\item $\rm{o}(||u||_{\xHone([0,t])}) =
\rm{o}(||u||_{\xHone([0,T+\delta])}) = \rm{o}(a_1,\ldots,a_n)$,
\item For $i=1,\ldots,n-1$~:
$dE_i^t(0)=dE_i^T(0)+\rm{O}(t-T)$.
\end{itemize}
Therefore we get, for $i=1,\ldots,n-1$~:
$$dE_i^t(0).u=a_i+\rm{o}(a_1,\ldots,a_n,t-T)$$
and thus~:
\begin{equation*}
\begin{split}
E_t(u)\ =\ (t+a_1+\rm{o}(a_1,\ldots,a_n,t-T),\ 
&a_2+\rm{o}(a_1,\ldots,a_n,t-T),\ \ldots,\ \\
&a_{n-1}+\rm{o}(a_1,\ldots,a_n,t-T),\  E_n^t(u))
\end{split}
\end{equation*}

Let us first solve the system of equations~:
$$E_1^t(u)=T,\ E_2^t(u)=\cdots=E_{n-1}^t(u)=0$$
We get from the \it{Implicit Function Theorem}~:
$$\forall i\in\{1,\ldots,n-1\}\quad a_i=f_i(T-t,a_n)$$
and moreover~:
$$a_1=T-t+\rm{o}(T-t,a_n)\quad \rm{ and for }i=2,\ldots,n-1~:
a_i=\rm{o}(T-t,a_n)$$

Let us now calculate the last coordinate~:
$$E_n^t(u)\ =\
E_n^t(0)+dE_n^t(0).u+\inv{2}d^2E_n^t(0).(u,u)+\rm{o}
(||u||^2_{\xHone}) \ =\ 
Q_t(u)+\rm{o}((T-t)^2,a_n^2)$$
Moreover~:
\begin{equation*}
Q_t(u)=Q_t(a_1(e_1^T+v)+a_2e_2^T+\cdots+a_{n-1}e_{n-1}^T+a_nw)
\end{equation*}
Hence
\begin{equation*}
\begin{split}
E_n^t(u)
&=(T-t)^2Q_t(e_1^T+v)+a_n^2Q_t(w)+a_n(T-t)\check{Q}_t(e_1^T+v,w)
+\rm{o}((T-t)^2,a_n^2)\\
&=g(t,a_n)
\end{split}
\end{equation*}
where $\check{Q}$ denotes the symmetric bilinear form associated
to the quadratic form $Q$.
On the one part~: $g(T,a_n)=a_n^2Q_T(w)+\rm{o}(a_n^2)>0$ if $a_n$
is small enough (by definition of $w$).
On the other part~: $g(t,0)=(T-t)^2Q_t(e_1^T+v)+\rm{o}((T-t)^2)$.
From the definition of $v$, we have~: $Q_T(e_1^T+v)<0$. Hence if
$t$ is close enough to $T$, but not equal to $T$, we get~:
$g(t,0)<0$.
Now applying the \it{Mean Value Theorem} to $g$ we can assert
that there exist $t$ close to T and $a_n$ small such that
$g(t,a_n)=0$, that is $E_n^t(u)=0$.
Actually there exist four such couples $(t,a_n)$
with $t>0$ or $t<0$, and $a_n>0$ or $a_n<0$.

\begin{rmrk} \label{rem34,6}
For such a couple $(t,a_n)$
we can say more about the asymptotics of $t$ with respect to
$a_n$, which will be useful in the sequel. To
this aim let us solve the following equation, where $a_n$ is
fixed~:
\begin{equation} \label{eqg}
g(t,a_n)=0
\end{equation}
We have~:
\begin{itemize}
\item $Q_t(e_1^T+v)=Q_T(e_1^T+v)+\rm{O}(T-t)$,
\item $Q_t(w)=Q_T(w)+\rm{O}(T-t)$,
\item
$\check{Q}_t(e_1^T+v,w)=\check{Q}_T(e_1^T+v,w)+\rm{O}(T-t)$.
\end{itemize}
Hence equation (\ref{eqg}) becomes~:
$$Q_T(e_1^T+v)(T-t)^2+a_n\check{Q}_T(e_1^T+v,w)(T-t) +a_n^2Q_T(w)
+\rm{o}((T-t)^2,a_n^2)=0$$
If $a_n$ is fixed, close to $0$, we easily find two solutions
$t_1<T$ and $t_2>T$ such that both satisfy~:
$$T-t=\beta a_n+\rm{o}(a_n),\ \beta\in\xR\backslash\{0\}$$
\end{rmrk}

To end the proof we have to study the cost $C_t(u)$ and
prove that $t$ and $a_n$ can be chosen so that this cost is lower
than $C_T(0)$. Recall vectors $v,w\in\xR e_0^T\oplus F$ previously
defined so that~:
\begin{equation} \label{ineqvw}
Q_T(e_1^T+v)<0\quad\rm{and}\quad Q_T(w)>0
\end{equation}
We can write~:
$$v=\lambda e_0^T+v_1,\ w=\mu e_0^T+w_1\quad\rm{ where }
v_1,w_1\in F$$
Up to adding $\varepsilon e_0^T$, $\varepsilon $ small, we can
assume that $\mu\neq 0$. This
does not change anything in inequalities (\ref{ineqvw}), nor in
our previous reasoning. Let us denote $C(t,u)=C_t(u)$. We can
expand~:
$$C_t(u)=C_T(0)+\alpha(t-T)
+\f{\partial C}{\partial u}(T,0).u +\rm{o}(T-t,||u||)$$
where $\alpha=\f{\partial C}{\partial t}(T,0)$. Moreover by
definition of $u$ and $e_i^T$~:
$$\f{\partial C}{\partial u}(T,0).u\ =\ 
dC_T(0).u\ =\ a_1\lambda+a_n\mu$$
On the one part~: $a_1=T-t+\rm{o}(T-t,a_n)$. On the other part
from Remark \ref{rem34,6} there exists $\beta\neq 0$ such that
$a_n=\beta(t-T)+\rm{o}(t-T)$. Hence we obtain~:
$$C_t(u)=C_T(0)+(t-T)(\alpha-\lambda+\beta\mu)+\rm{o}(t-T)$$
Now among the four solutions of Remark \ref{rem34,6}, we can
choose $t$ and $a_n$ such that $\alpha-\lambda+\beta\mu\neq 0$
and $(t-T)(\alpha-\lambda+\beta\mu)<0$. In these conditions~:
$$C_t(u)<C_T(0)$$
which ends the proof.
\end{proof}


\paragraph{{\it{2. Final time fixed}}}
\begin{dfntn}
\begin{itemize}
\item
The trajectory $\gamma$ is said to be \it{$\xCzero$-cost-minimal} on
$[0,T]$ if there exists a $\xCzero$-neighborhood of $\gamma$ such
that for any trajectory $q$ contained in this neighborhood, with
$q(0)=\gamma(0)$ and $q(T)=\gamma(T)$, we have~: $C_T(v)\geq
C_T(0)$, where $v$ is the control associated to $q$.
\item
The trajectory $\gamma$ is said to be
\it{$\xLinfty$-cost-minimal} on $[0,T]$
if there exists a neighborhood of $0$ in $\xLinfty([0,T])$ such
that, for any trajectory $q$ associated to a control $v$ of this
neighborhood, with $q(0)=\gamma(0)$ and $q(T)=\gamma(T)$, we
have~:
$C_T(v)\geq C_T(0)$.
\end{itemize}
\end{dfntn}

\no We have the following~:

\begin{thrm} \label{coutmintff}
The trajectory $\gamma$ is
$\xCzero$-cost-minimal if and only if $T<t_c$. Moreover, $\gamma$
is not $\xLinfty$-cost-minimal if $T>t_c$
(whereas $\gamma $ is $\xCzero$-time-minimal if and only if
$T<t_{cc}$), where $t_{cc}$ and $t_c$ denote the two types of
first conjugate times of $\gamma$.
\end{thrm}

Hence in this case, the times at which $\gamma$ ceases
to be minimizing are different in the time-optimal problem and
cost-optimal problem~: $\gamma$
ceases to be $\xCzero$-time-optimal \it{before} it ceases to be
$\xCzero$-cost-optimal (since $t_{cc}<t_{c}$,
see Lemma \ref{tempsconj}).

\begin{proof}
The proof is quite similar to proof of Theorem \ref{coutmintfnf}
(but simpler)
and is only sketched.
In the coordinates $(x_1,\ldots,x_n)$ of Lemma
\ref{formenormaleexcept}, we write $E_T=(E_1^T,\ldots,E_n^T)$,
where $d^2E^T_n(0)=0$,
and $dE_1^T(0),\ldots$, $dE_{n-1}^T(0),dC_T(0)$ are independant.
We decompose~:
$$\xLtwo([0,T])=\xR e_1\oplus\cdots\oplus\xR e_{n-1}\oplus\xR e_0\oplus F$$
where $dE_i^T(0).e_i=1$ and $dC_T(0).e_0=1$. By definition of
$t_c$~:
\begin{itemize}
\item if $T<t_c$~: $Q_{/ \xR e_0\oplus F}$ is positive definite,
\item if $T>t_c$~: $Q_{/ \xR e_0\oplus F}$ is indefinite.
\end{itemize}
where $Q=d^2E_n^T(0)$.
Hence if $T>t_c$ there exist $v,w\in \xR e_0\oplus F$ such that
$Q(v)>0$ and $Q(w)<0$. Up to adding $\varepsilon e_0$ we can
assume that projections of $v,w$ on $\xR e_0$ are not trivial.
Consider the control~:
$$u=a_1e_1+\cdots+a_{n-1}e_{n-1}+\lambda v+\mu w$$
and try to solve equations~: $E_i(u)=\delta_i^1T, i=1\ldots n$.

First from the Implicit Function Theorem we get, solving the
$n-1$ first coordinates (i.e. the controllable part of the
system)~:
$$a_i=f_i(\lambda,\mu)=\rm{o}(\lambda, \mu), i=1\ldots n-1$$
Now the last coordinate is~:
\begin{equation*}
\begin{split}
E_n(u)&=Q(u)+\rm{o}(||u||^2_{\xLtwo})\\
&=\lambda^2Q(v)+\mu^2Q(w)+\lambda\mu\check{Q}(v,w)+\rm{o}(
\lambda^2, \mu^2)
\end{split}
\end{equation*}
Using the Mean Value Theorem we obtain~:
\begin{itemize}
\item $\exists \lambda_1<0,\mu_1>0\ /\ E_n(
a_1e_1+\cdots+a_{n-1}e_{n-1}+\lambda_1 v+\mu_1
w)=0$
\item $\exists \lambda_2>0,\mu_2>0\ /\ E_n(
a_1e_1+\cdots+a_{n-1}e_{n-1}+\lambda_2 v+\mu_2
w)=0$
\end{itemize}
There exists $i=1 \rm{ or }2$ such that $dC_T(0).(\lambda_i
v+\mu_i
w)\neq 0$, and therefore there exists a trajectory steering $0$
to $\gamma(T)$ in time $T$ with a cost strictly lower than
$C_T(0)$.
\end{proof}


\subsection{Optimality of abnormal trajectories for
sub-Riemannian systems of rank 2}
With the notations of Section \ref{sectionsrcase} we have the
following result, which can be proved in the same way as Theorem
\ref{bordSR}.

\begin{thrm} \label{thrmoptimaliteSR}
Under assumptions of Theorem \ref{bordSR},
the abnormal reference trajectory $\gamma$ is $\xCzero$-optimal
for the sub-Riemannian system (\ref{systemesousR}) if and only if
it is $\xCzero$-time-minimal for its associated affine system
(\ref{affi}). Moreover $\gamma$ is abnormal for this
affine system~; actually $\gamma$ is $\xCzero$-optimal
if $T<t_{cc}$ and is not $\xLinfty$-optimal if $T>t_{cc}$.
\end{thrm}

In particular \it{conjugate times are the same along
$\gamma$ for the sub-Riemannian system (\ref{systemesousR}) and
its associated affine system (\ref{affi}).}
Therefore the whole formalism
that was introduced for affine systems
(the differential operators $D_1, D_2$) is still valid in
sub-Riemannian geometry. Hence the conjugate time of the
sub-Riemannian problem can be computed using an algorithm. This
result makes a link between works of \cite{BK} and \cite{AS},
\cite{AS1}.

\begin{xmpl}
The Martinet case (see Section \ref{secteursSR})
is in dimension 3, hence $t_{cc}=+\infty$
(see Remark \ref{remtcc2}).
The abnormal trajectory is optimal on $\xR^+$.
\end{xmpl}

\begin{rmrk}
As proved in \cite{Agra} the $\xCzero$-optimality is in
sub-Riemannian geometry equivalent to the optimality in the
$\xLtwo$ topology on controls.
\end{rmrk}

\begin{rmrk}
If $T$ is small enough (depending on the choice of the Riemannian
structure, and lower than $t_{cc}$),
then as first noted by \cite{AS} $\gamma$ is moreover
\it{globally} optimal among all sub-Riemannian trajectories
steering $0$ to $\gamma(T)$.
\end{rmrk}

\begin{rmrk}
It should be noted that the loss of optimality holds in
$\xLinfty$. Hence using the definitions introduced in Section
\ref{secteursSR}, the $\xLtwo$-sector plays no role in the
optimality of the reference abnormal trajectory. The loss of
optimality holds in the $\xLinfty$-sector.
\end{rmrk}


\begin{acknowledgement}
I would like to thank my teacher B.
Bonnard for many relevant ideas and advices.
\end{acknowledgement}

\end{document}